\documentclass{article}

\usepackage{amssymb}
\usepackage{amsmath}
\usepackage{graphics}

\let\<\langle
\let\>\rangle
\font\bigcmsy=cmsy10.pk scaled 2000
\def\bigtimes{\mathop{\,\vrule width0pt depth2pt height8pt
            \smash{\lower2pt\hbox{\bigcmsy\char'002}}\,}\limits}

\begin{document}

\begin{center}
\Large{\textbf{Quasipotency and $g$-potency in free constructions.}}
\end{center}
\begin{center}

\textbf{Vladimir V. Yedynak}

\end{center}

\begin{abstract}

In this paper we investigate the properties of quasipotency and $g$-potency for free products with amalgamation.


\end{abstract}

\section{Introduction.}

Definition 1. A group $G$ is called quasipotent if for for each element $g$ of $G$ there exists a natural number $k$ such that for each natural number $n$ there exists a finite index subgroup $H_n$ such that $\<g\>\cap H_n=\<g^{kn}\>$.

The property of quasipotency was investigated in [1] where different sufficient conditions were found for quasipotency of free products, HNN extensions and extensions. Quasipotency was investigated in [1] for to obtain the classes of cyclic subgroup separable groups.

Definition 2. Consider a group $G$ with a fixed element $g$. We say that $G$ is $g$-potent if for each natural number $n$ there exists a finite index subgroup $H_n$ such that $\<g\>\cap\<H_n\>=\<g^n\>$.

If the property from the definition 2 is true for each element of a group $G$ then $G$ is called potent. In [2] it was proved that free groups are potent. In [5] there were investigated links between $g$-potency and residual finiteness of free product with cyclic amalgamation.

In [4] it was derived the property of \textsl{regular quotient} at a pair $\{U, V\}$ of a group $G$ reminding the property of quasipotency. A group $G$ is said to have regular quotient at a pair $U, V$ of its elements if there exists number $r$ so that for each $n$ there exists a homomorphism $\varphi$ of $G$ onto a finite group such that the orders of $\varphi(U)$ and $\varphi(V)$ are equal to $rn$. Later this property for free groups was generalized in [3] where it was proved that free groups are omnipotent. A group $G$ is called omnipotent if for each $n$ and for each elements $g_1,..., g_n$ of $G$ such that $g_i$ and $g_j, i\neq j$ have no conjugate nonunit powers there exists a number $k$ such for each ordered sequence of natural numbers $l_1,..., l_n$ there exists a homomorphism of $G$ onto a finite group such that the order of $\varphi(u_i)$ equals $kl_i$. This definition generalizes the property of quasipotency.

In this work we prove the following theorems.

\textbf{ Theorem 1. }\textsl{ Let $A$ be a quasipotent group and $B$ is a subgroup of $A$ such that $A$ is $B$-separable that is $B$ is the intersection of finite index subgroups of $A$. Consider the finite set $I$ and isomorphisms $\varphi_i: A\rightarrow A_i, i\in I$. Put $B_i=\varphi_i(B)$. Then the group $G={\ast\atop B_j=\varphi_{jk}(B_k)}A_i,\ i,\ j,\ k\in I,\ j\neq k, \varphi_{jk}=\varphi_j\varphi_k^{-1}$ is quasipotent.}

\textbf{ Theorem 2. }\textsl{ Let $A_i, i\in I$ be residually finite groups and $N$ is the kernel of the natural homomorphism of the group $G={\ast\atop i\in I}A_i$ onto a direct product of groups $A_i, i\in I$. Then for each element $h$ of $N$ the group $G$ is $h$-potent. }

\textbf{ Theorem 3. }\textsl{ Let $A_i, B_j, i\in I, j\in J$ be a residually finite groups. Consider elements $u\in
F={\ast\atop i\in I}A_i, v\in G={\ast\atop j\in J}B_j$ not belonging to subgroups which are conjugate to free factors and isomorphism $\varphi: \<u\>\rightarrow\<v\>, \varphi: u\mapsto v$. Then the group $H=F{\ast\atop \<u\>\simeq^{\varphi}\<v\>}G$ is residually finite. }

\section{Auxiliary notations.}

Consider the group $G$ generated by the set $\{g_i| i\in I\}$ and acting on the set $X$. Construct the graph $\Gamma$ such that the set of vertices of $\Gamma$ coincides with $X$. If $p\circ g_i=q$ and $p\circ g_i^{-1}=r$ for some $p, q, r\in X$ and for some $i\in I$ then the vertices $p$ and $q$ are connected by the edge going from $p$ into $q$ laballed by $g_i$ and $p$ and $r$ are connected by the edge going from $r$ into $p$ labelled by $g_i$. We consider that for each edge $f$ of $\Gamma$ with a label there exists the inverse edge $f'$ without a label such that the begin point of $f$ coincides with the end of $f'$ and the beginning of $f'$ is the end of $f$. The graph $\Gamma$ is oriented and edges with labales are positively oriented edges of $\Gamma$ and the rest edges are negatively oriented edges of $\Gamma$.

Definition 1. The graph $\Gamma$ which is constructed as above is called the action graph of the group $G$ with respect to its generating set $\{g_i| i\in I\}$.

We shall omit further the references about the generating set if it is fixed. Besides throughout the paper considering the action graphs of free product we shall always suppose that the set of generators coincides with the union of free factors.

Consider the action graph $\Gamma$ of the group $G=\<g_i| i\in I\>$.

The label of the edge $e\in\Gamma$ is denoted as Lab $(e)$. The label of the path $S=e_1...e_n$ of $\Gamma$ is the element of the group $G$ equal to $\prod_{i=1}^n$ Lab $(e_i)'$ where Lab $(e_i)'$ equals Lab $(e_i)$ in case $e_i$ is positively oriented or Lab $(e_i')^{-1}$ otherwise where $e_i'$ is the edge inverse to $e_i$. The beginning and the end of the edge $e$ will be denoted as $\alpha(e)$ and $\omega(e)$ correspondingly.

Definition 2. Consider the element $u$ of the group $G$ which does not equal to unit. Let $\Gamma$ be the action graph of the group $G=\<g_i|i\in I\>$. Fix the vertex $p$ in $\Gamma$. Then the \textsl{$u$-cycle $T$} of the graph $\Gamma$ going from $p$ is the set of paths $\{S_i=\{e_j^i| j\in J_i\}|i\in J\}$ satisfying the following conditions:

1) $\alpha(S_l)=p$;

2) there exists a one to one correspondence between the paths $S_l$ and the presentations of $u$ in generators $\{g_i|i\in I\}: u=g_{i_1}^{\varepsilon_1}...g_{i_k}^{\varepsilon_k}, $
min$(|i_j-i_{j+1}|, \varepsilon_j+\varepsilon_{j+1})>0,
\varepsilon_i\in\{-1; 1\}, 1\leqslant j\leqslant k, 1\leqslant
i\leqslant k$. Besides if $\varepsilon_j=1$ then edge $e_{kn+j}^l$ is positively oriented and its label equals $g_{i_j}$ for each natural $n$; if $\varepsilon_j=-1$ then $e_{kn+j}^l$ is negatively oriented and the label of the edge which s inverse to $e_{i_j}$ equals $g_{i_j}$ for each natural $n$; in case the path $S_l$ is finite and is composed of $r$ edges we consider that indices are modulo $r$;

3) there is no closed subpath $K$ in the path $S_l$ which differs from $S_l$ and satisfies the conditions 1), 2);

4) the path $S_l$ is the maximal path on the entry which satisfies the conditions 1)-3);

The paths which compose a $u$-cycle will be called as representatives of a $u$-cycle.

Definition 3. Suppose that some representative of the $u$-cycle $T$ is composed from the finite set of edges and has the label $u^k$. Then we shall say that the \textsl{length} of the $u$-cycle $T$ is equal to $k$.

Suppose we have a group $G$ which is defined by its presentation $G=<x_i, i\in I| R_j, j\in J>$. It is obvious that a graph $\Gamma$ is the action graph of the group $G=\<x_i| i\in I\>$ if and only if it satisfies the following conditions:

1) for each vertex $p$ of $\Gamma$ and for each $i\in I$ there exist exactly one edge with label $x_i$ going from $p$ and exactly one edge going into $p$ labelled by $x_i$;

2) for each $j\in J$ and for each vertex $p$ of $\Gamma$ the $R_j$-cycle of $\Gamma$ going from $p$ has length 1;

Using this remark we can describe action graphs of the group $G={\ast\atop B_k=B_j} A_i$ corresponding to actions of $G$ on $X$ which are so that $A_i$ act freely on $X$ where $i\in I, B_i<A_i$ and $B_j$ is isomorphic to $B_k$ for each $j, k\in I$. The action graph $\Gamma$ of such kind satisfies the following properties.

(1) for each vertex $p$ of $\Gamma$ and for each $c\in {\bigcup\atop i\in I}A_i$ there exists exactly one edge with label $c$ going into this vertex and there exists exactly one edge labelled by $c$ going from $p$;

(2) for each vertex $p$ of $\Gamma$ and for each $i\in I$ the maximal connected subgraph $A_j(p)$ whose positively oriented edges are labelled by elements from $A_i$ is the Cayley graph of the group $A_i$ with the generating set $\{A_i\}$; let $B$ be the amalgamated subgroup of $G$ then the subgraph $B(p)$ is defined analogically as $A_i(p)$;

In what follows we shall use the following notation. Consider the graph $\Gamma$ and its subset $T$. Consider $t$ copies of $\Gamma: \Gamma_1,..., \Gamma_t$. Then $T^i$ is the set of $\Gamma_i$ which corresponds to $T$ in $\Gamma$.

Let $\Gamma$ be the graph with properties (1), (2). Fix the set of vertices $S\in\Gamma$ and the function $f: S\rightarrow I$. Besides for each $p$ and $q$ from $S$ at least one of the following conditions is true: $f(p)\neq f(q), A_{f(p)}(p)\neq A_{f(q)}(q)$. Consider $t$ copies of $\Gamma: \Delta_1,..., \Delta_t$ from which we shall construct the new action graph. For each $p\in S$ and for each edge $e$ whose internal point belongs to $B(p)$ and whose label belongs to $A_{f(p)}\setminus B$ we delete edges $e^i$ for all $i=1,...,t$. If the edge $e$ goes into the vertex $p'\in B(p)$ from the vertex $q$ we connect the vertices $p'^i$ and $q^{i+1}$ by the edge going from $q^{i+1}$ whose label is equal to Lab($e^i$) and if the edge $e$ goes away from the vertex $p'$ into the vertex $q$ then we connect vertices $p'^i$ and $q^{i+1}$ by the edge which goes away from $p'^i$ whose label equals Lab($e^i$) (indices are modulo $t$). Put $S=\{s_1,...., s_n\}$. We obtained the new graph which satisfies the properties (1), (2) so it is the action graph ant it will be referred to as $\gamma_t(\Gamma; s_1,..., s_n;f(s_1),..., f(s_n))$.

Having a group $G={\ast\atop B_k=B_j} A_i$ we shall call the graphs with properties (1), (2) as free action graphs of the group $G$.

Consider the free action graph $\Gamma$ of $G={\ast\atop B_k=B_j} A_i$. Fix some element $k$ from $I$. Let $S=e_1...e_n$ be a path of $\Gamma$ passing through the vertex $p$. Let $l^+_{p, k}(S)$ be the number of edges in $S$ going into the vertices of the subgraph $B(p)$ and having labels from $A_k\setminus B^{\ast}$ and $l^-_{p, k}(S)$ is the number of edges of $S$ which goes away from vertices of $B(p)$ and which have labels from $A_k\setminus B^{\ast}$ (some of them perhaps are equal since the set $\{e_1,..., e_n\}$ of edges may contain repeatable edges). Put $l_{p, k}(S)=|l^+_{p, k}(S)-l^-_{p, k}(S)|$.

Definition 3. Consider the vertex $p$ of $\Gamma$ and $k\in I$ and the representative $T'$ of the $u$-cycle $T$ of finite length. If $l_{p, k}(T')\neq0$ then we shall say that $T$ is proper relative to $p$ and $k$.

In this paper we use the symbol $(m, n)$ to designate the greatest common divisor of numbers $m$ and $n$.

Definition 4. Consider two vertices $p$ and $q$ of the graph $\Gamma$ which belong to one component of $\Gamma$. Then the
\textsl{distance} between $p$ and $q$ is the minimal length of the finite path which connects $p$ and $q$. The length of the finite path is the number of edges which compose this path.

Definition 5. Consider the graph $\Gamma$, the cycle $S=e_1...e_n$ of $\Gamma$ and the nonnegative integer number $l$. We say that $S$ does not have \textsl{$l$-near vertices} if for each $i, j, 1\leqslant i<j\leqslant n$ the distance between the vertices
$\alpha(e_i), \alpha(e_j)$ is greater or equal than min$(|i-j|, n-|i-j|, l+1)$.

\section{Auxiliary lemmas.}

Throughout this section we consider that we have a fixed group $G={\ast\atop B^{\ast}=B_i} A_i$.

\textbf{ Lemma 1.}\textsl{ Let $\Gamma$ be the free action graph of the group $G$. Fix the vertex $p$ of $\Gamma$ and the element $k\in I$ and let $u$ be the cyclically reduced element not belonging to any free factor. Consider the $u$-cycle $S$ of finite length in $\Gamma$ going from the vertex $q$ and its representative $S'$ passing through the vertex $p$. Then for each $t$ the $u$-cycle $T$ going from the vertex $q_1$ of the graph $\Delta=\gamma_t(\Gamma; p; k)$ has length equal to $t/(t, l_{p, k}(S))l(S)$ where $l(S)$ is the length of the $u$-cycle $S$. }

\textsl{ Proof. }

Let $\Delta_1,...,\Delta_t$ be the copies of the graph $\Gamma$ from which the graph $\gamma_t(\Gamma; p; k)$ is composed of. Fix the edge $e$ of the graph $\Gamma$ which either has the label different from an arbitrary element of $A_k\setminus B^{\ast}$ or does not contain the vertex $p$ in the capacity of the terminal vertex.

Notice that there exists the singular composition of the representative $S'$ as the union of paths

$$
(i)\ S=S_1\bigcup...\bigcup S_k,
$$

where each path $S_i$ has the following presentation $S_i=e_{1i}...e_{n_i, i}$ where $\alpha(e_{1i}), \omega(e_{n_i, i})\in B(p)$ and $\alpha(e_{j, i})\notin B(p), \omega(e_{j, i})\notin B(p)$ for each $j, 1<j<n_i$.

Consider the path $S_i=f_1...f_l$ where $\alpha(f_1)=p', \omega(f_l)=p''\in B(p)$. Put $S_i'^j$ be the path in $\Delta$ with the presentation $S_i'^j=rf_2^j...f_{l-1}^js$ where $r$ is either the $f_1^j$ if Lab$(f_1)\notin A_k$ or it is the new edge which was appended during the construction of the graph $\Delta$ instead of the edge $f_1^j$ in case Lab$(f_1)\in A_k$. The edge $s$ is defined analogically.

We shall work with the representative $T'$ of the $u$-cycle $T$ going from the vertex $q^1$ such that $T'$ corresponds to the same notation of $u$ as $S'$.

Changing $u$ on its appropriate cyclic permutation we may consider that the first edge of $S'$ goes away from the vertex of $p$ and the cycle $S'$ traverses the path $S_1$ first next $S_2$ etc. and eventually returns to the vertex $p$. Some of the paths $S_i$ are perhaps traversed by $S'$ several times.

Suppose that during the pass-by of the representative $T'$ we traverse the path $S_i'^j$ and $\alpha(S_i'^j)=p'^r$ then $\omega(S_i'^j)=p''^s$ for some natural $r, s$. Besides if $f_1$ and $f_l$ do not belong to $A_k(p)$ then $r=s=j$ if both edges $f_1, f_l$ belong to $A_k(p)$ then $r=s=(j-1)($mod\ $t)$. In case $f_1\in A_k(p)$ and $f_l\notin A_k(p)$ then $(r+1)($mod\ $t)=j=s$ and finally if $f_1\notin A_k(p)$ while $f_l\in A_k(p)$ then $r=j=s+1($mod\ $t)$. Hence we deduce that making the pass-by of $T'$ from the vertex $p_i$ and passing the path whose label equals the label of the path $S$ we come to the vertex from the subgraph $B(p_{i+\varepsilon})$ where $\varepsilon=l_{p, k}(S)\cdot$sgn$(l_{p, k}^+(S)-l_{p, k}^-(S))$ (indices at  $p$ are modulo $t$). There exists the minimal on modulus natural number $q$ such that $t | q\varepsilon$. If $t=b(t,
l_{p, k}(S)), l_{p, k}(S)=c(t, l_{p, k}(S))$, then $q=b=t/(t, l_{p, k}(S))$. Notice that $q$ is equal to the ration of the length of $u$-cycle going from the vertex $p_i$ to length of the $u$-cycle $S$. Lemma 1 is proved.

\textbf{ Lemma 2.}\textsl{ Consider the free action graph $\Gamma$ of the group $G$ so that for each vertices $p, q$ of $\Gamma$ the subgraphs $A_i(p), A_j(q), i\neq j$ have no two mutual vertices which belong to different subgraphs $B(r), B(r')$. Fix the cyclically reduced element $u\in G\setminus\{\cup_{i\in I}A_i\}$. Suppose that for each vertex $p$ of $\Gamma$ and for each $k\in I$ $l_{p, k}(S)=0$ where $S$ is an arbitrary representative of an arbitrary $u$-cycle. Consider that the length of each $u$-cycle of $\Gamma$ is equal to the same number. Fix the representative $S=e_1...e_n$ of some $u$-cycle of $\Gamma$ and put $\alpha(e_1)=p$, Lab$(e_2)\in A_k$, Lab$(e_{n-1})\in A_s$. Then in the graph $\Delta=\gamma_2(\Gamma; \alpha(e_2), \alpha(e_n); k, s)$ the length of each $u$-cycle equals the length of a $u$-cycle from $\Gamma$. Besides the vertex $p^1$ is the beginning of the representative of the $u$-cycle which is proper relative to $p^1$ and $k$.}

\textsl{ Proof. }

Suppose that Lab$(e_1)\in A_t,$ Lab$(e_n)\in A_r$. Notice that the condition of the theorem involves that the subgraphs $B(\alpha(e_2)), B(\alpha(e_n))$ are different. Because otherwise $A_t(\alpha(e_1))$ and $A_r(\alpha(e_n))$ possesses two common vertices $\alpha(e_1)$ and $\alpha(e_n)$ not belonging to the subgraph $B(\alpha(e_n))$.

Let $\Delta_1, \Delta_2$ be two copies of the graph $\Gamma$ from which the graph $\Delta$ is composed of.


Let $T$ be the representative of the $u$-cycle going from $p^1$, $T$ corresponds to the same notation of $u$ as $S$. If Lab$(e_1)\in A_t$ then there exists exactly one edge in $T$ going from the vertex of the subgraph $B(p^1)$ and belonging to $A_t(p^1)$. This edge is $e_1^1$. Hence $l_{p^1,t}(T)=1$.

Let's prove now that the length of each $u$-cycle in $\Delta$ coincides with the length of the $u$-cycle $S$. Let $T$ be a representative of the $u$-cycle in $\Gamma$ going from the vertex $r$ and passing through the vertex from $B(\alpha(e_2^i))$ or $B(\alpha(e_n^i))$ for some $i=1, 2$.

Put $\alpha(e_2)=b, \alpha(e_n)=c$. Consider the unique factorization of $T$ as the union of paths:

$$
(i)\ T=T_1\bigcup...\bigcup T_m,
$$
where each path $T_i$ has the presentation $T_i=e_{1i}...e_{n_i, i}$ such that $\alpha(e_{1i}), \omega(e_{n_i, i})\in B(b), \alpha(e_{j, i})\notin B(b), \omega(e_{j, i})\notin B(b), 1<j<n_i$. Notice that $l_{c, s}^{\pm}(T)=\sum_{i=1}^ml_{c, s}^{\pm}(T_i)$. As follows from the proof of lemma 1 going from the vertex $b_i$ and traversing the path $T_j$ we come to the vertex from the subgraph $B(p_{i+w})$ where $w=(l^+_{b, k}(T_j)-l^-_{b, k}(T_j)+l^+_{c, s}(T_j)-l^-_{c,
s}(T_j))($mod $2)$. Since $l_{c, s}^{\pm}(T)=\sum_{i=1}^ml_{c, s}^{\pm}(T_i), l_{b, k}^{\pm}(T)=\sum_{i=1}^ml_{b, k}^{\pm}(T_i)$ and $T$ is not proper relative to any vertex and element of the set $I$ the $u$-cycle of $\Delta$ going from $r_1$ has length which is equal to the length of the $u$-cycle containing $S$. Lemma 2 is proved.

\textbf{ Lemma 3.}\textsl{ Consider natural numbers $n, l_1,..., l_k$. Let $d=(l_1,..., l_k)$ be the greatest common divisor of numbers $l_1,..., l_k$. Then the least common multiple $F$ of numbers $nd/(nd, l_1),..., nd/(nd, l_k)$ equals $n$. }

\textsl{ Proof. }

Let $p$ be the prime number such that for some natural number $m$ it is true that $p^m\ |\ n, p^{m+1}\ \nmid\ n$. Then $n=p^mr, (r, p)=1$. Since it is true that $d=(l_1,..., l_k)$ then there exists number $i$ and there exists nonnegative integer $c$ such that $l_i=p^cq, d=p^cs, (q, p)=(s, p)=1, q=st$. Then $nd/(nd, l_i)=p^mr\cdot p^cs/(p^mr\cdot p^cs, p^cst)=p^mr/(p^mr, t)=p^m\cdot r/(r, t)$ because $t=q/s$ and $p$ are coprime. For each $j\neq i$ it true that $l_j=p^{c+\varepsilon}st', \varepsilon\geqslant0, (t',p)=1$ and $nd/(nd, l_j)=p^mr/(p^mr,
p^\varepsilon t')$. Hence $p^{m+1}\nmid nd/(nd, l_j)$. Therefore $p^m\mid F, p^{m+1}\nmid F$ which means that $n\mid F, (F/n, n)=1$. For each $i$\ $d\mid(nd, l_i)$ hence $F\mid n$. Lemma 3 is proved.

\textbf{ Lemma 4. }\textsl{ Let $A_i, i\in I$ be a finite set of finite groups such that each group $A_i$ contains the subgroup $B_i$ and all the subgroups $B_i$ are isomorphic: $B_j=\varphi_{jk}(B_k), G={\ast\atop B_j=\varphi_{jk}(B_k)}A_i, i, j, k\in I, j\neq k$. Consider some cyclically reduced element $u$ from $G$ which does not belong to any free factor. Consider also the finite free action graph $\Gamma$ of $G$ so that the lengths of all $u$-cycles of $\Gamma$ are equal to $m$. Besides for each vertices $p, q$ of $\Gamma$ the subgraphs $A_i(p), A_j(q)$ have no two mutual vertices which do not belong to one subgraph $B(r)$ for each distinct $i$ and $j$. Then for each natural number $n$ there exists the homomorphism of $G$ onto a finite group such that the order of $u$'s image equals $mn$. }

\textsl{Proof. }

Suppose there exists the vertex $p$ of $\Gamma$ and the element $k\in I$ so that there exists the representative of a $u$-cycle proper relative to $p$ and $k$. Let $S_1,..., S_l$ be all representatives of $u$-cycles proper relative to $p$ and $k$. Put $d=(l_{p, k}(S_1),..., l_{p,
k}(S_l))$. Then according to lemma 1 in the graph $\Delta=\gamma_{nd}(\Gamma; p; k)$ the length of each $u$-cycle equals one of the numbers: $m, m\cdot nd/(nd, l_{p, k}(S_1)),..., m\cdot nd/(nd, l_{p, k}(S_l))$. According to lemma 3 the least common multiple of these numbers equals $mn$. Since the group $G$ acts on the set of vertices of $\Delta$ there exists the homomorphism $\varphi: G\rightarrow S_w$, where $w$ is the number of vertices in $\Delta$. The value of $\varphi(u)$ is the product of independent cycles and the least common multiple of their lengths is equal to $mn$ that is $|\varphi(u)|=mn$.

Suppose that for each vertex $p$ of $\Gamma$ and for each element $k\in I$ each representative of each $u$-cycle is not proper relative to $p$ and $k$. Consider the $u$-cycle $S=e_1...e_k$ and $\alpha(e_1)=p$, Lab$(e_1)\in A_t$, Lab$(e_2)\in A_k$, Lab$(e_{n-1})\in A_s$. According to lemma 2 in the graph $\Delta=\gamma_2(\Gamma; \alpha(e_2), \alpha(e_n); k, s)$ there exists the $u$-cycle whose representative is proper relative to $p^1, t$ and the length of each $u$-cycle in $\Delta$ equals $m$. Thus the proof can be finished due to the first case.

\section{Proof of theorems.}

\textsl{ Proof of theorem 1. }

Consider the element $u\in G$. We may consider that $u$ is cyclically reduced and has infinite order. Suppose that $u\in A_s$ for some $s$. Since the group $A_s$ is quasipotent there exists the natural number $k_u$ such that for each natural number $n$ there exists the homomorphism $\phi_s$ of $A_s$ onto a finite group $K_s$ such that the image of $u$ has order $k_un$. For each $j\in
I\setminus\{s\}$ define the homomorphism $\phi_j:A_j\rightarrow K_j$, $\phi_j=\phi_s\varphi_s\varphi_j^{-1}$. Put $L_i=\phi_i(B_i)$. Then there exists the homomorphism $\phi:G\rightarrow K={\ast\atop L_j=\phi_{jk}(L_k)}K_i, i, j, k\in I, j\neq k, \phi_{jk}=\phi_j\varphi_{jk}\phi_k^{-1}$ which is defined by the following condition. If $v=a_1...a_n\in G$ is a reduced notation of the element $v$, $a_j\in A_{i_j}\setminus B_{i_j}$ then $\phi(v)=\phi_1(a_1)...\phi_n(a_n)$, if $a_i\in
A_k$ then $\phi_i(a_i)=\phi_k(a_i)$. Besides the element $\phi(u)$ has order $k_un$. Since the group $K$ is residually finite there exists the homomorphism of $K$ onto a finite group such that the image of $\phi(u)$ has order $k_un$.

Consider the case when $u\notin{\bigcup\atop i\in I}A_i$. We fix some reduced notation for the element $u$: $u=u_1...u_n, u_j\in A_{i_j}\setminus B_{i_j}$. Since the group $A_{i_j}$ is residually finite with respect to the entry in the subgroup $B_{i_j}$ there exists the homomorphism of $A_{i_j}$ onto a finite group such that the image of each element from $\{u_1,..., u_n\}\cap A_{i_j}$ does not belong to the image of the subgroup $B_{i_j}$. Let $N_{i_j}$ be the kernel of this homomorphism. Put $M_{i_j}=\varphi_{ij}^{-1}(N_{i_j}), M={\bigcap\atop j}M_{i_j}$. Since for each $j=1,...,n$ $\mid A:M_{i_j}\mid<\infty$, it is true that $\mid A:M\mid<\infty$. Put $N_i=\varphi_i(M),L_i=\varphi_i(B_i), i\in I$. Then we may define the natural homomorphisms $\phi_i:A_i\rightarrow K_i=A_i/N_i, i\in I$. These homomorphisms can be completed till the homomorphism $\phi:G\rightarrow K={\ast\atop L_j=\phi_{jk}(L_k)}K_i, i, j, k\in I, j\neq k, \phi_{jk}=\phi_j\varphi_{jk}\phi_k^{-1}$ be the following way. If $v=a_1...a_i\in G$ is a reduced notation of the element $v$, $a_j\in A_{i_j}\setminus B_{i_j}$ then $\phi(v)=\phi_{i_1}(a_1)...\phi_{i_n}(a_n)$, if $a_j\in A_k$ then $\phi_{i_j}(a_i)=\phi_k(a_i)$. Thus we may consider that the free factors $A_i$ are finite. Then the group $G$ is residually finite. Consider the homomorphism $\varphi: G\rightarrow K$ where $K$ is a finite group and $\varphi$ is so that $({\bigcup\atop i\in I}A_i\cup(\cup_{i, j}(\{A_i\}\times\{A_j\})))\bigcap$Ker$\varphi=\{1\}$. Then it is possible to consider that the subgroups $A_i$ are embedded into the group $K$ for each  $i$. Let $\Gamma$ be the Cayley graph of the group $\varphi(G)$ with the generating set $\{{\bigcup\atop i\in I}A_i\}$. Put $k_u=\ \mid\varphi(u)\mid$. According to lemma 4 for each natural number $n$ there exists the homomorphism of $G$ onto a finite group such that the image of $u$ has order $k_un$. It means that the group $G$ is quasipotent. Theorem 1 is proved.

\textsl{ Proof of theorem 2. }

Consider the natural homomorphism $\varphi: G\rightarrow {\times\atop i\in I}A_i$, $N=\ $ker$\varphi$. Notice that we may consider that the set $I$ is finite. Let $u=u_1...u_k$ be the nonunit element of the subgroup $N$. Using the residual finiteness of the groups $A_i$ for each $i\in I$ we define the homomorphism $\varphi_i$ of $A_i$ onto a finite group such that if  $u_s\in A_i$ then $\varphi_i(u_s)\neq 1$. So we may consider that free factors are finite. Let $\Gamma$ be the Cayley graph of the group ${\times\atop i\in I}A_i$ with the generating set $\{{\cup\atop i\in I}A_i\}$. Then $\Gamma$ is a finite free action graph of $G$ in which all $u$-cycles have length 1. Besides according to the definition of the direct product it follows that for each vertex $p\in\Gamma$ and for each $i, j\in I, i\neq j$ the subgraphs $A_i(p), A_j(p)$ have the singular common vertex. According to lemma 4 for each natural number $n$ there exists the homomorphism of $G$ onto a finite group such that the image of $u$ has order $n$. Theorem 2 is proved.

The proof of the next theorem illustrates the application of the notion of $g$-potency to the investigation of residual finiteness of the free product with cyclic amalgamation.

\textsl{ Proof of theorem 3. }

Notice that we may consider that the sets $I$ and $J$ are finite. Consider the element $w$ from $H$. For each $i\in I, j\in J$ define $\alpha_i, i\in I$ as the set of elements from $A_i$ which composes the normal forms of elements $u, w$. Let $\beta_j, j\in J$ be the set of elements of $B_j$ which composes the normal forms of $v, w$. Since $A_i, B_j$ are residually finite there exist homomorphisms $\varphi_i, \psi_j$ from $A_i, B_j$ onto finite groups such that $\alpha_i\cap$ ker $\varphi_i, \beta_j\cap$ ker $\psi_j$ are empty. These homomorphisms can be completed till the homomorphisms $\varphi_F: F\rightarrow \ast_{i\in I}\varphi_i(A_i), \varphi_G: G\rightarrow
\ast_{j\in J}\psi_j(B_j)$. Since elements $u'=\varphi_F(u), v'=\varphi_G(v)$ have infinite orders there exists the homomorphism $\varphi_H: H\rightarrow (\ast_{i\in I}\varphi_i(A_i)){\ast\atop \<u'\>\simeq^{\phi}\<v'\>}(\ast_{j\in J}\psi_j(B_j)), \phi=\varphi_G\varphi\varphi_F^{-1}$. Thus we may consider that all groups $A_i, B_j$ are finite for each $i, j$. It means that there exists the natural number $n$ such that $u^n, v^n$ belong to the Cartesian subgroups of groups $F$ and $G$ correspondingly. So according to the theorem 2 the groups $F$ and $G$ are $u^n$- and $v^n$-potent correspondingly. Since it is known that the cyclic subgroup separability is inherited by free products the groups $F$ and $G$ are cyclic subgroup separable. It is sufficient just to refer to the following theorem that was proved in [5].

\textbf{Theorem}: \textsl{ The group $A{\ast\atop\<u\>_{\infty}\simeq^{\varphi}\<v\>_{\infty}}B, \varphi(u)=v$ is residually finite if the groups $A, B$ are residually finite and $A$ is $u^n$-potent and $u$-separable and $B$ is also $v^n$-potent and $v$-separable for some natural $n$.}

Theorem 3 is proved.

\textbf{Theorem 4.}\textsl{ Consider the finitely generated group $G$ which contains the subgroups $K, F$ such that $G=KF$, $K$ is finite and $F$ is free. Then $G$ is quasipotent.}

\textsl{Proof.}

Since $G$ is residually finite it is sufficient to consider the elements of infinite order. Consider the element $u\in G\setminus\{\cup_{g\in G}g^{-1}Kg\}$. Fix the basis in $F: x_1,..., x_m$. Since $G$ is virtually free $G$ is subgroup separable. So there exists the homomorphism $\varphi$ of $G$ onto a finite group such that $\varphi_{K}$ is injective and in the graph $\Gamma=Cay(\varphi(G); \{\varphi(K)\}\cup\{\varphi(x_1),..., \varphi(x_m)\})$ the representatives of $u$-cycles have no 1-near vertices. Put $k=|\varphi(u)|$. For each natural number $n>1$ we construct the action graph of the group $G=\<K, x_1,..., x_m\>$ --- $\Gamma_n$ by the following way. Fix in $\Gamma$ an arbitrary vertex $p$. Consider $n$ copies of the graph $\Gamma: \Delta_1,..., \Delta_n$. Let $R$ be the subgraph in $\Gamma$ containing $p$ which is the Cayley graph of the group $K$ with the set of generators $\{K\}$. For each $i, j, 1\leqslant i\leqslant n, 1\leqslant j\leqslant m$ and for each vertex $r$ from $R^i$ consider the edge $q^i$ with label $x_j$ going into the vertex $r^i$. Delete the edge $q^i$ and connect the vertex $\alpha(q^i)$ with the vertex $r^{i+1}$ (indices are modulo $n$) by the edge $f_i$ with label $x_j$, which goes away from the vertex $\alpha(q^i)$. The obtained graph is the action graph of the group $G=\<K, x_1,..., x_n\>$. Since all representatives of $u$-cycles in $\Gamma$ have no 1-near vertices then the constructible graph contains a $u$-cycle of length $kn$. Theorem 4 is proved.

\textbf{Theorem 5.} \textsl{Consider the group $G=A\ast B$ and $u=ab\in G$ where $A$ and $B$ are finite and $a\in A, b\in B, |a|=|b|>1$. Then $G$ is $u$-potent.}

\textsl{Proof.}

Consider the natural number $n$. The idea of the proof is to construct the free action graph of the group $G$ in which each $u$-cycle has length which is either equal to $n$ or to 1. Consider the Cayley graph $P=Cay(A; \{A\})$ of the group $A$ and the Cayley graph $Q=Cay(B; \{B\})$ of the group $B$. Consider $n$ copies of $P$ and $Q: P_1,..., P_n, Q_1,..., Q_n$. Put $m=|a|=|b|$. Consider the cycles $\Lambda_a$ in $P$ and $\Lambda_b$ in $Q$ whose labels are equal to $a^m, b^m$ correspondingly and whose edges have labels $a$ and $b$ correspondingly. Let $p$ and $q$ be the nonnegative natural numbers such that $m=2+p+q$. Put $\Lambda_a=e_1...e_m, \Lambda_b=f_1...f_m$. In what follows all upper indices are modulo $n$. We consider that $\omega(e_1^i)=\alpha(f_1^i), \omega(f_1^i)=\alpha(e_1^{i+1})$ for each $i=1,..., n$. Besides $\omega(e_{p+2}^i)=\alpha(f_{p+2}^i), \omega(f_{p+2}^i)=\alpha(e_{p+2}^{i+1})$ for all $i$. Furthermore $\omega(e_{j+1}^i)=\alpha(f_{m-j}^i), \omega(f_{j+1}^i)=\alpha(e_{m-j}^{i+1})$ when $i=1,..., n$ and $j=1,..., p$. Notice that due to the embedding of $A\ast B$ into $(A\times B)\ast(A\times B)$ we may consider that $|A|=|B|$ and hence the number of cycles of $P$ with label $a^m$ whose edges have label $a$ coincides with the number of cycles in $Q$ that have label $b^m$ and whose edges have label $b$. If we do the similar procedure for graphs $P_1,..., P_n, Q_1,..., Q_n$ for each pair of cycles with labels $a^m$ and $b^m$ we obtain the action graph of the group $G$ which is so that each $u$-cycles has length which equals either $n$ or 1. Theorem 5 is proved.

\begin{figure}[h]
\centering
\includegraphics*[9cm,75 mm]{1-st.bmp}
\end{figure}

\newpage
\begin{center}
\large{Acknowledgements.}
\end{center}

The author thanks A. A. Klyachko for valuable comments and discussions.

\begin{center}
\large{References.}
\end{center}

1. \textsl{J. Burillo, A. Martino.} Quasi-potency and cyclic
subgroup separability. Journal of Algebra, Volume 298, Issue 1,
Pages 188-207

2. \textsl{P. Stebe.} Conjugacy separability of certain free
products with amalgamation. Trans. Amer. Math. Soc. 156 (1971).
119--129.

3. \textsl{Wise, Daniel T.} Subgroup separability of graphs of free groups with cyclic edge groups.
Q. J. Math. 51, No.1, 107-129 (2000). [ISSN 0033-5606; ISSN 1464-3847]

4. \textsl{Graham A. Niblo.} H.N.N. extensions of a free group by \textbf{z} which are subgroup separable. Proc. London Math. Soc. (3), 61(1):18--23, 1990.

5. \textsl{Allenby R. B. J. T., Tang C. Y.} The residual finiteness
of some one-relator groups with torsion. // J. Algebra 71 (1981),
\textbf{1}. 132--140

\end{document}